\RequirePackage{fix-cm}
\documentclass[smallextended]{svjour3}       
\smartqed  
\usepackage{graphicx}
\usepackage{amsfonts, amsmath, float, supertabular, array}
\usepackage[english]{babel}
\usepackage{blindtext}
\usepackage[ruled,linesnumbered]{algorithm2e}

\makeatletter

\newcommand{\Rmnum}[1]{\expandafter\@slowromancap\romannumeral #1@}
\makeatother
\newcommand{\rad}{{\rm rad}}
\renewcommand{\mid}{{\rm mid}}

\begin{document}

\title{An efficient algorithm for global interval solution of nonlinear algebraic equations and its GPGPU implementation}


\titlerunning{An efficient algorithm for global interval solution of nonlinear algebraic equations}        

\author{LIN Dang         \and
        CHEN Liangyu 
}


\institute{LIN Dang \at
              Shanghai Key Lab of Trustworthy Computing \\
              \email{51151500028@stu.ecnu.edu.cn}           
           \and
           CHEN Liangyu \at
              Shanghai Key Lab of Trustworthy Computing\\
              \email{lychen@sei.ecnu.edu.cn}
}

\date{Received: date / Accepted: date}

\maketitle

\begin{abstract}
Solving nonlinear algebraic equations is a classic mathematics problem, and common in scientific researches and engineering applications. There are many numeric, symbolic and numeric-symbolic methods of solving (real) solutions. Unlucky, these methods are constrained by some factors, e.g., high complexity, slow serial calculation, and the notorious intermediate expression expansion. Especially when the count of variables is larger than six, the efficiency is decreasing drastically. In this paper, according to the property of physical world, we pay attention to nonlinear algebraic equations whose variables are in fixed constraints, and get meaningful real solutions. Combining with parallelism of GPGPU, we present an efficient algorithm, by searching the solution space globally and solving the nonlinear algebraic equations with real interval solutions. Furthermore, we realize the Hansen-Sengupta method on GPGPU. The experiments show that our method can solve many nonlinear algebraic equations, and the results are accurate and more efficient compared to traditional serial methods.
\end{abstract}

\section{Introduction}
\label{intro}
Construction theory and algorithm for algebraic equations is very classical and important. Many problems in scientific researches and engineering fields, eventually are transferred into equations to be solved. For linear equations, there are lots of mature researches and tools. By the common desktop computers, one can easily solve linear equations with thousands of variables or higher. For nonlinear algebraic equations, there are also many symbolic and numeric methods of solving solutions. The methods of exact solutions for nonlinear algebraic equations can be divided into two categories, symbolic algorithm and symbolic-numerical algorithm. The symbolic methods include Groebner base~\cite{14}, Wu's method~\cite{24}, resultant elimination method~\cite{25} and others. In recent years, some symbolic-numerical methods have also been developed, to preserve the accuracy of symbol calculation and enjoy the benefits of numerical computation. Many of the above methods are also implemented into specific software packages, such as Groebner package in Maple, WSolve~\cite{21}, Discoverer~\cite{20,26}, GAS~\cite{25}, Epsilon~\cite{23} and more. However, with the rapid development of physical world, such as applications in space control, intelligent transportation, program verification, internet of things and etc, it has brought big challenges for solving the nonlinear algebraic equations. The current methods and tools have the following limitations. The first is insufficient computation. Most of the specific packages described above, currently, can only be executed with one core. Worse, it is not easy to extend these packages for execution on cluster computers or other high performance environments. The second is insufficient memory utilization caused by intermediate expression expansion. In successive computation, it generates many intermediate symbolic expressions, which cannot be omitted due to preserve strict accuracy. There is also no floating number truncation and rounding error. Thus, these intermediate expressions make the physical memory exhausted quickly and computation aborted.

Currently, the computer hardware has passed the multi-core age, and stridden to the many-core age. Many desktop computers have equipped with multiple CPU cores, usually four or eight cores, and the Intel CPU with MIC architecture has 32 cores or more. More lucky, CPU is not the sole calculation part in modern computers. GPU (Graphic Processing Unit) Computing has also been widely applied and popularized. The optimized GPGPU program can work very fast. However, most of traditional programs in symbol calculation still follow the old serial executive rule even with good CPU and GPU devices. Therefore, combining parallel computing and symbolic computation, and using parallel computing to accelerate symbolic computation process, have become more important in current researches, e.g., Wang's Parallel polynomial operations on SMPs~\cite{19}, Lin's parallel computation for polynomial GCD~\cite{22,17}, Sato's on the parallel computation of comprehensive groebner systems~\cite{16}, Moreno Maza's multithreaded parallel implementation of arithmetic operations modulo a triangular set~\cite{17}, Kobayashi's work of the parallel implementation of the formal verification language-Isabelle. Using GPGPU technology to solve symbolic computation problems, is also emerging with the utilization of thousands of cores. Some researchers have make progress on theoretical analysis and calculation under many-core environments, like~\cite{15,18}.

In this paper, based on the branch and bound method~\cite{2} and interval arithmetic~\cite{1,10,28}, we present an efficient algorithm to get potential real  solutions for nonlinear algebraic equation with integer or floating number coefficients in GPGPU environment. It's easily observe that, many problems in physical world can be transformed into nonlinear algebraic equations with two characteristics. One is each variable has limited value ranges. The other is only real solutions are emphasized and required. Therefore, we apply the simple idea of branch and bound method, and do successively iterative partition in the global solution space to remove unfeasible parts. We also use the Hansen-Sengputa method~\cite{4,5,7,13} to accelerate the convergence of solution space. The main advantage of our algorithm is complete and sound, namely, it can quickly find all subdivisions who potentially contain real solutions. We globally search the solution space, and use reliable interval arithmetics without any floating number truncation or rounding error to keep interval solutions right. Moreover, whole calculations are executed not only in numerical way but also in parallel, with the support from thousands of cores in GPGPU. This also avoids the notorious problem of intermediate expression expansion.

The structure of the paper is listed as follows. In section 2, we introduce the fundamental knowledge of interval arithmetics and Hansen-Sengputa method. In section 3, the main algorithm and related sub-algorithms are presented. In section 4 we present experiments and data statistics. The conclusion is in section 5.

\section{Preliminaries}
\label{basic}
\subsection{Polynomial equation}
A polynomial is an expression consists of variables (or indeterminates) and coefficients, involving only the operations of addition, subtraction, multiplication, and variables with non-negative integer exponents. A typical example of a single indeterminate $x$ is $x^{2}-3x+4$.

A polynomial function is a function that can be defined by evaluating a polynomial. For example, the function $P$
\begin{equation}\label{eq:1}
P=a_{0}+a_{1}X+\cdots +a_{n}X^{n}.
\end{equation}

The corresponding polynomial function can construct the equation
\begin{equation}\label{eq:2}
f_{P}(x)=a_{0}+a_{1}x+\cdots +a_{n}x^{n}=0.\
\end{equation}

If $r \in \mathbb{R}$, and $f_{P}(r)=0$, it can call $r$ the root of the polynomial equation, or the zero of the associated function. The relation of the root of polynomial function and the polynomial is: if $r\in \mathbb{R}$ is a root of $f_{P}$ if and only if $X-r$ divides $f_{P}$, that is, there exist another polynomial $Q(x)$, such as: $f_{P}=(X-r)Q(x)$. If a positive integer $k$, which $f_{P}=(X-r)^{k}Q(x)$, then call $r$ a multiple root of $f_{P}$.

A real root means the root of an equation is a real number. Then the real root isolation of the polynomial system is presenting the intersection intervals of the real number line, which contain all the real roots of the polynomial, and each interval has one and only one real root.
\subsection{Interval arithmetic}
 A real interval is a set of real numbers with the property that any number that lies between two numbers is also included in the set. For example, let $X$ be a real interval and let
\begin{equation}\label{eq:3}
X=[\underline{x}, \overline{x}]=\{x \in \mathbb{R}: \underline{x}\leq x\leq \overline{x} \},\
\end{equation}
where $\underline{x}$ is the lower bound, $\overline{x}$ is the upper bound respectively. The set of all interval numbers can also be defined as
\begin{equation}\label{eq:4}
\mathbb{IR}=\{[\underline{x}, \overline{x}]: \underline{x}, \overline{x}\in \mathbb{R}, \underline{x}\leq \overline{x} \}.
\end{equation}

Let $c$ be the midpoint of an interval, we have \begin{equation}\label{E:2.5}
c=\mid(x)=(\underline{x}+\overline{x})/2.
\end{equation}

Let $r$ be the radius of an interval, we also have\begin{equation}\label{E:2.6}
r=\rad(x)=(\overline{x}-\underline{x})/2.
\end{equation}

The basic arithmetic operations of interval computing are defined in $\mathbb{IR}$, so the result of interval calculation contains all possible real solutions, $x=[\underline{x}, \overline{x}], y=[\underline{y}, \overline{y}]$. These basic operations are listed as follow.
\begin{equation}\label{eq:5}
x+y=[\underline{x}+\underline{y}, \overline{x}+\overline{y}],\
\end{equation}
\begin{equation}\label{eq:6}
x-y=[\underline{x}-\overline{y}, \overline{x}-\underline{y}],\
\end{equation}
\begin{equation}\label{eq:7}
\begin{split}
x\times y=[\min\{\underline{x}\underline{y}, \overline{x}\underline{y}, \underline{x}\overline{y}, \overline{x}\overline{y}\},
  \max\{\underline{x}\underline{y}, \overline{x}\underline{y}, \underline{x}\overline{y}, \overline{x}\overline{y}\}],\
\end{split}
\end{equation}
\begin{equation}\label{eq:8}
1/x=[1/ \overline{x}, 1/ \underline{x}] \quad \mbox{if}\ \underline{x} > 0 \ \mbox{or}\ \overline{x}<0,\
\end{equation}
\begin{equation}\label{eq:9}
x \div y = x \times 1/y.
\end{equation}

In the above division operations, the case of zero division is undefined. We make use of extended interval arithmetic introduced by Hanson~\cite{8} and Kahan~\cite{29} and deal the division with \

\begin{eqnarray}\label{eq:10}
x/y = \left\{ \begin{array}{ll}
 [\overline{x}/\underline{y}, +\infty]  & \textrm{: $  \overline{x}  \leq 0  \  \& \  \overline{y}=0 $}, \\

[-\infty, \overline{x}/\overline{y}] \bigcup [\overline{x}/\underline{y}, +\infty]  & \textrm{: $ \overline{x}  \leq 0  \  \& \  \underline{y} < 0 < \overline{y}$}, \\

 [-\infty, \overline{x}/\overline{y}]  & \textrm{: $  \overline{x}  \leq 0 \  \& \ \underline{y}=0 $}, \\

  [-\infty, +\infty]  & \textrm{: $  \underline{x}  < 0  < \overline{x}$}, \\

 [-\infty, \underline{x}/\underline{y}]  & \textrm{: $  \underline{x}  \geq 0  \  \& \  \overline{y}=0 $}, \\

 [-\infty, \underline{x}/\underline{y}] \bigcup [\underline{x}/\overline{y}, +\infty]  & \textrm{: $ \underline{x}  \geq 0  \  \& \  \underline{y} < 0 < \overline{y}$}, \\

  [\underline{x}/\overline{y}, +\infty]  & \textrm{: $  \underline{x}  \leq 0  \  \& \  \underline{y}=0 $}.\\

\end{array} \right.
\end{eqnarray}

\definition{Let $f$ be an arithmetic representation in $\mathbb{R}[x_1, x_2, \ldots, x_n]$, and transform all operands to intervals, thus all arithmetic operations can be seen as the corresponding interval operations, denote as $F$. Then $F: I(\mathbb{R})^n \rightarrow I(\mathbb{R})$ can be called an interval evaluation, simply, $Y \supseteq {f(x):x \in X}$ define the evaluation of $f(x)$ in $X$.}

Generally, the methods of interval evaluation of a function in a given region can be classified into three categories:

1) Using the algebra method or analysis directly, calculate the interval value of the function in a given region.

2) Dividing an given interval into a number of subintervals, calculate the interval function value on each subinterval, and choose the maximum and minimum as the result of the original function evaluation bound.

3) Using the optimization methods for function in a given range, and choose the maximum and minimum as evaluation bound.

\definition{Let $f$ be an arithmetic representation in $\mathbb{R}[x_1, x_2, \ldots, x_n]$. Given a space $I = [[\underline{x_1}, \overline{x_1}], [\underline{x_2}, \overline{x_2}], \ldots, [\underline{x_n}, \overline{x_n}]]$} with $n$ tuples, the function $f$ can be calculated by interval evaluation and the result is also in interval, namely, $[\underline{f},\overline{f}]$. If $0 \in [\underline{f},\overline{f}]$, it means $I$ may contain real solutions, and we call the tuple $C$ as an configuration of $f$; otherwise, $I$ is an infeasible part of $f$.

\definition{An interval matrix is a matrix whose elements are intervals. We denote the interval matrix $[X]$ with its elements $[X]_{ij}$, and have
\begin{equation}\label{eq:11}
[X]=\begin{bmatrix}
[\underline{x_{11}}, \overline{x_{11}}] &[\underline{x_{12}}, \overline{x_{12}}]&\cdots & [\underline{x_{1n}}, \overline{x_{1n}}]\\
[\underline{x_{21}}, \overline{x_{21}}] &[\underline{x_{22}}, \overline{x_{22}}] &\cdots & [\underline{x_{2n}}, \overline{x_{2n}}]\\
 \vdots   & \vdots & \ddots  & \vdots  \\
[\underline{x_{n1}}, \overline{x_{n1}}] &[\underline{x_{n2}}, \overline{x_{n2}}] &\cdots & [\underline{x_{nn}}, \overline{x_{nn}}]
\end{bmatrix}.
\end{equation}}
\subsection{Hansen-Sengupta method}
Hansen-Sengupta method is an algorithm using interval arithmetic to compute and bound the zeros of nonlinear algebraic equations. It provides guaranteed bounds on all zeros in a given region. Hansen-Sengupta method is closely related to the interval Newton method~\cite{5,8,9}.\

Let $f$ be a function of $n$ variables, $x= (x_1,x_2,\ldots,x_n)^\mathrm{T}$, and function
\begin{equation}\label{eq:12}
f(x) = 0.
\end{equation}

We expand $f$ by using Taylor's theorem and expanding $f(y)$ about $x$, then obtain
\begin{equation}\label{eq:13}
f(x)+J(\xi)(z-x)=f(y)=0,
\end{equation}
where $J(\xi)$ is a Jacobian matrix evaluated at a point $\xi \in X$, and the point set $z$ contains all solution of $y$.
We would like to know the set $z$, but it's difficult to represent this set.

Moore~\cite{10} first found out to calculate the inverse of the Jacobian matrix $J(X)$ in interval Newton method. But it require an interval matrix $M$ containing every interval element. Hansen~\cite{14} pointed out that it is not necessary to find an interval inverse, instead, a Gaussian elimination procedure could be used as in the real counterpart.

Define $J_c$ as the center of $J(X)$, and each element of $J_c$ is the midpoint of the corresponding interval element of $J(X)$.
A nonlinear algebraic equations with interval coefficients such as (\ref{eq:13}) is best solved by multiplying an approximate inverse of $J_c$. Let $B$ be this approximation matrix.

We thus rewrite (\ref{eq:13}) as
\begin{equation}\label{eq:14}
Bf(x)+BJ(X)(z-x)=0.
\end{equation}

The products $Bf(x)$ and $BJ(X)$ are computed in interval arithmetic to bound rounding errors.

Krawczyk introduced a variation of interval Newton method to avoid Gaussian elimination on interval matrix by not attempting to obtain a sharp solution of (\ref{eq:14}). Based on this method, we can compute the box with
\begin{equation}\label{eq:15}
K(X)=x-Bf(x)+[I-BJ(X)](X-x).
\end{equation}

If a solution $y$ is contained in a box $X$, then it is also contained in $K(X)$. Since $K(X)$ may not be contained in $X$, we can use the iteration
\begin{equation}\label{eq:16}
X^{(i+1)}=X^{i} \cap K(X^{i}) (i=0, 1, \ldots, n).
\end{equation}

The Krawczyk method computes the box $K(X)$ in each iteration. However, it is not the smallest box.

We now present Hansen-Sengupta method. The box is generally smaller than $K(X)$. Each iteration of Hansen-Sengupta method tends to produce a greater reduction of the current box than Krawczyk's method does, also fewer steps are required for convergence.

Let $g=Bf(x)$ and $P=BJ(X)$, thus (\ref{eq:15}) is rewritten as
\begin{equation}\label{eq:17}
g+P(z-x)=0.\end{equation}

Note that the interval matrix $P$ is
\begin{equation}\label{eq:18}
P=L+D+U,
\end{equation}
where the matrices $L$, $D$, and $U$ are lower triangular, diagonal, and upper triangular, respectively. We thus rewrite (\ref{eq:17}) as
\begin{equation}\label{eq:19}
Y=x-D^{-1}[g+L(X'-x)+U(X-x)], 
\end{equation}

\begin{equation}\label{eq:20}
X'=Y\cap X.
\end{equation}
After each component $Y_i(i=1, 2, \ldots, n, )$ is obtained, it intersects with $X_i$ so that the new result $X^{\prime}_i=Y_i \cap X_i$ can be used in obtain $Y_{i+1}, Y_{i+2}, \ldots, Y_n$, thus we compute componentwise, for $i=1,2,\ldots, n$,

\begin{equation}\label{eq:21}
\begin{split}
Y_i=x_i-(D_{ii})^{-1}[g_i+ \sum_{j=1}^{i-1}P_{ij}(X'_j-x_j)+ \sum_{j=i+1}^nP_{ij}(X_j-x_j)],
 \end{split}
\end{equation}

\begin{equation}\label{eq:22}
X'_i=Y_i \cap X_i.
\end{equation}

This step is done for each $i = 1,2,\ldots,n$ and the process is iterated until the new box is sufficiently small.

Note that even though $P$ is supposed to approximate the identity matrix, the interval $D_{ii}$ may contain zero. We simply use extended interval arithmetic to compute $Y_i$. The intersection (\ref{eq:22}) then produces a finite result.
\subsection{Gauss-Jordan algorithm}
Gauss-Jordan Algorithm~\cite{6} is a classic method for matrix inversion. Let $I_n$ be an identity matrix of size $n$. Consider a matrix $A$ of size $n$, and an augment $A$ to get the matrix $[C]=[A|I_n]$,
thus we can do column rank elimination of Gauss-Jordan method on $C$ for $k=1, 2, \ldots, n-1$.

Firstly, select the main element by column maximum$|A_{ik}|$. Exchange the elements of  $k$-th row and $ik$-th row of $[C]=[A|I_n]$.

Secondly, calculate the main element
\begin{equation}\label{eq:23}
C_{kj} \leftarrow C_{kj}/A_{kk}, (j=k,k+1 \ldots, 2n).\
\end{equation}

Finally, do elimination calculation
\begin{equation}\label{eq:24}
\begin{split}
C_{ij} \leftarrow C_{ij}-A_{ik} \times C_{kj}, (i=1, 2, \ldots, n; i\not=k; j=k+1,k+2, \ldots, 2n).
\end{split}
\end{equation}

By doing elimination calculation on each row successively, the left half becomes the $I_n$, and the right half becomes the desired inverse of A.

\section{A nonlinear algebraic equations interval algorithm based on branch and bound method}
\label{algorithm}
\subsection{The framework of our algorithm}
 Our goal is to solve the problem of nonlinear algebraic equations whose each variable is in a given region. Assume the variables are $x_1, x_2, \ldots, x_n$, and their regions are $[\underline{d_i}, \overline{d_i}]$ respectively. So an initial solution space can be constructed easily by combining these regions together. According to the branch and bound method, we can divide the solution space to $2^n$ rectangular subdivisions~\cite{27} by splitting the region to two half parts for each variable and making combinations on these new parts. Each of these subdivisions can be deemed as a potential solution space. Through certain mathematical calculation, some subdivisions can be proved no roots and removed. The next step is using the remained subdivisions as original inputs. We can do the division and determination process iteratively until the interval width of subdivisions is sufficiently small.

 Since the computation on intervals is "loose", the final result may have too many configurations to be understood. Moreover, the intervals are very small, so that many subdivisions are nearly overlapped. Thus, we can do backtracking on these configurations to get more readable result.

\begin{figure}[H]
  \includegraphics[width=1\textwidth]{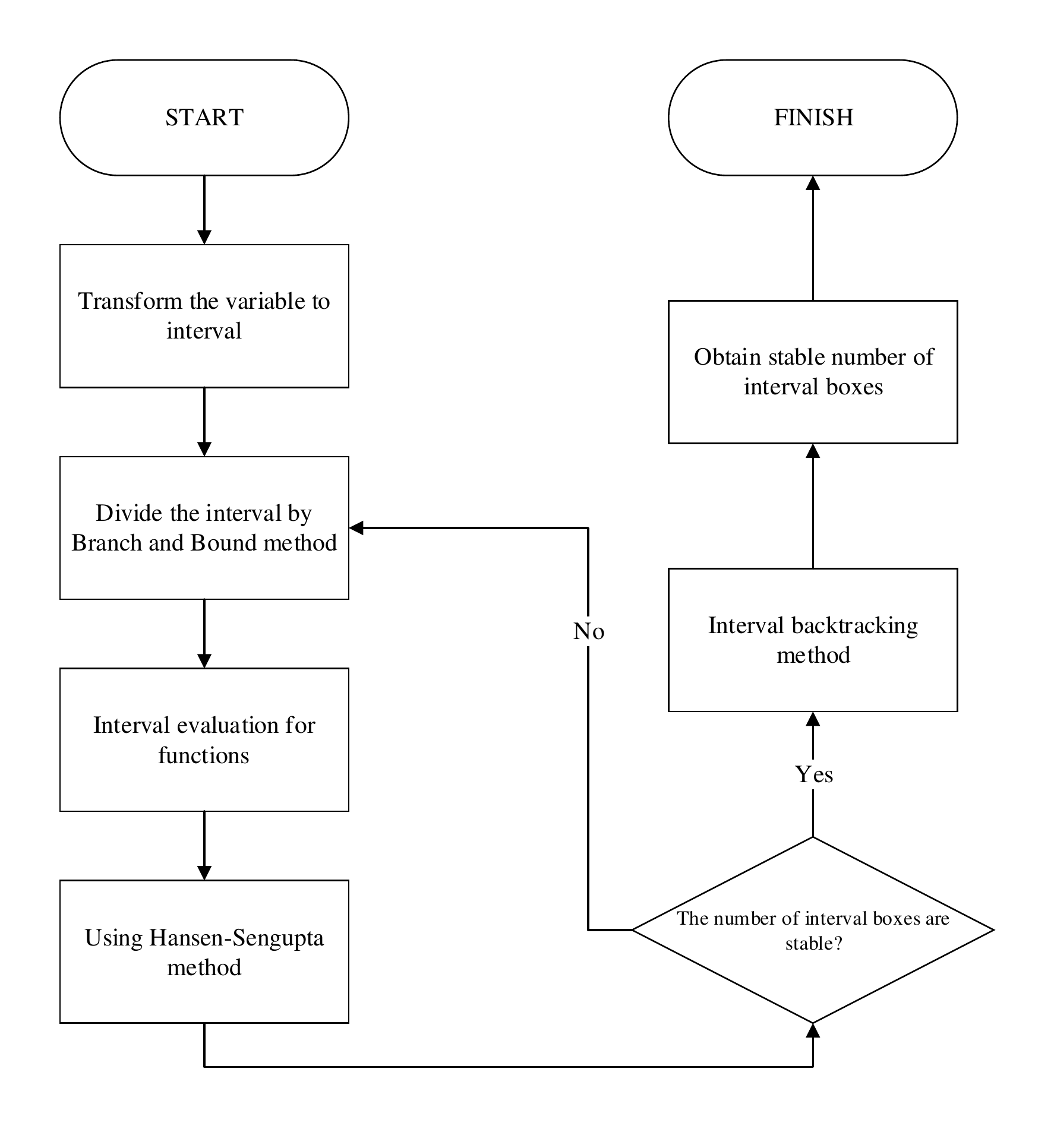}
\caption{The overall flowchart of Algorithm 1}
\label{fig1}       
\end{figure}

Our algorithm is described as Algorithm~\ref{alg1} and its flowchart is showed in Fig.\ref{fig1}.

\begin{algorithm}[H]
\caption{Solving nonlinear algebraic equations system based on branch and bound method}
\label{alg1}
\KwIn{The algebraic equations $f(x_1,x_2,\ldots,x_n)$, its variable regions and the desired width or result intervals}
\KwOut{The feasible solution set $M_k$, where each element is an configuration with $n$ intervals and may contain feasible solutions.}
According to the original regions, construct an interval configuration as the initial solution set $M_0$.

\textbf{For} each element from $M_0$ \textbf{Do}\linebreak
 (\Rmnum{1}) divide the configuration into $2^n$ subdivisions by using branch and bound method;\linebreak
 (\Rmnum{2}) make evaluation on subdivisions and remove the infeasible ones;\linebreak
 (\Rmnum{3}) use Hansen-Sengupta methods to compute and bound the zeros of nonlinear functions, and remove the infeasible ones;\linebreak
 (\Rmnum{4}) add the remained subdivisions into $M_1$.

Do the operations described in the step 2 iteratively and get the $M_0, M_1, \ldots, M_k$,
calculation is terminated if one of the conditions is satisfied:\linebreak
(\Rmnum{1}) {the $M_k$ is an empty set, then we can assert that there is no real solution for $f(x_1,x_2,\ldots,x_n)$;\linebreak
(\Rmnum{2}) $M_k$ is satisfied with the desired interval width, then $M_k$ is the final result.}
\end{algorithm}

In successive iterations, for one configuration with $n$ intervals, it generates $2^n$ new subdivisions according to the branch and bound method. So we need apply the Hansen-Sengupta method to remove more infeasible branches. While using Hansen-Sengupta method, we involve several sub-algorithms including interval matrix multiplication, Gauss-Jordan method for matrix inverse.

We obtain feasible interval configurations by iterate algorithms~\ref{alg2},~\ref{alg3}, then merge the interval configurations by backtracking Algorithm~\ref{alg5} and isolate the final roots within intervals.

\subsection{Polynomial interval evaluation and branch and bound method on GPGPU}
Firstly, we consider the problem of polynomial interval evaluation. Our research goal is nonlinear algebraic equations  with integer or floating number coefficients. To ensure the result strict accurate, all calculations are done in interval style. Detailedly, on GPGPU, we need carefully deal with the upper and lower bounds of basic calculations. In \textit{CUDA} programming, there is two commands supported by IEEE standard to ensure the correctness of interval evaluation. One is \textit{rd}, round-down, and the other is \textit{ru}, round-up. Therefore, the basic calculations in GPGPU follow these rules:
\begin{equation}\label{eq:25}
x+y=[\mbox{rd}(\underline{x}+\underline{y}), \mbox{ru}(\overline{x}+\overline{y})],
\end{equation}
\begin{equation}\label{eq:26}
x-y=[\mbox{rd}(\underline{x}-\overline{y}), \mbox{ru}(\overline{x}-\underline{y})],
\end{equation}
\begin{equation}\label{eq:27}
\begin{split}
x \times y=[\mbox{rd}(\min\{\underline{x}\underline{y}, \overline{x}\underline{y}, \underline{x}\overline{y}, \overline{x}\overline{y}\}),
   \mbox{ru}(\max\{\underline{x}\underline{y}, \overline{x}\underline{y}, \underline{x}\overline{y}, \overline{x}\overline{y}\})].
\end{split}
\end{equation}

Note that calculation in Algorithm~\ref{alg1} only involves polynomial evaluation, including addition, substraction and multiplication, but no any division.

To solve a nonlinear algebraic equation, the branch and bound method is a good way to remove infeasible parts from global solution space. For a configuration with $n$ interval tuples, we present Algorithm~\ref{alg2} to demonstrate one round of branch and bound method.

\begin{algorithm}[H]
\caption{One round of branch and bound method on GPGPU}
\label{alg2}
\KwIn{$f(x_1,x_2,\ldots,x_n)$, and the initial configuration $I=[[\underline{x_1}, \overline{x_1}],[\underline{x_2}, \overline{x_2}], \cdots, [\underline{x_n}, \overline{x_n}]] $.}
\KwOut{Configurations contain potentially feasible solution.}
Split each variable interval into two parts, namely, for $I$ from $n$ intervals to $2n$ intervals, $I_{11}=[\underline{x_1}, \mid(\underline{x_1}+\overline{x_1})], I_{12}=[\mid(\underline{x_1}+\overline{x_1}), \overline{x_1}, ], \ldots, I_{n1}=[\underline{x_n}, \mid(\underline{x_n}+\overline{x_n})], I_{n2}=[\mid(\underline{x_n}+\overline{x_n}), \overline{x_n}]]$.
Make combination on these new intervals for $n$ variables, and get new $2^n$ configurations.

Put $2^n$ configurations into $2^n$ threads of GPGPU, and each thread calculates $f(x_1,x_2,\ldots,x_n)$ with one configuration respectively.

According to the function interval evaluations, remove infeasible configurations (not contain zero) and output the remains.
\end{algorithm}

It is remarked that if $2^n$ exceeds the limit of GPGPU, the task can be divided into several jobs, which can be submitted repeatedly.

After several rounds of branch and bound method, the infeasible parts of global solution space are removed, while those feasible configurations potentially contain solutions are left.

\subsection{Hansen-Sengupta method based on  GPGPU}
Hansen-Sengupta method is a good method to determine whether a function has roots in a given box. If yes, this method can generate a more tighter box for roots of input function. We present the Algorithm~\ref{alg3} to demonstrate this method implementation on GPGPU.

\begin{algorithm}[H]
\caption{Hansen-Sengupta method on GPGPU}
\label{alg3}
\KwIn{$F=[f_1, f_2, \ldots, f_n] \in \mathbb{Q}[x_1, x_2, \ldots, x_n]$, the initial box $B=X_1 \times X_2 \times \cdots \times X_n$.}
\KwOut{The new box $B^*$ where\\
\quad \quad \quad \ 1) $\forall (x_1, x_2, \ldots, x_n) \in B, $ if $F(x_1, x_2, \ldots, x_n)=0$,\\ \quad \quad \quad \ then $(x_1, x_2, \ldots, x_n) \in B^*$.\\
\quad \quad \quad \ 2) if $B^* \subseteq B$, there exists $(x_1, x_2, \ldots, x_n) \in B^*$, \\ \quad \quad \quad \ and $ F(x_1, x_2, \ldots, x_n)=0$.
   \\}
Assume the Jacobian matrix of $F$ is $J$. Substitute the intervals $X_1, X_2, \ldots, X_n$ into $J$, then get a Jacobian matrix $J(X)$ whose entries are intervals.

Calculate center of $J(X)$, then calculate the approximate inverse $A=Inv(\mid(J(X)))$ with Gauss-Jordan method.

Calculate $M=A*J(X), R(x)=A*F(x)$ by interval matrix multiplication.\\
\textbf{For} $i=1, 2, \ldots, n$, \textbf{Do}\linebreak
$P_i=R(x)-\sum_{j=1}^{i-1}M_{ij}(X'_j-x_j)-\sum_{j=i+1}^{n}M_{ij}(X_j-x_j)$, \linebreak
$N_i(x, X)=x_i+P_i/M_{ii}$.

\textbf{For} $i=1, 2, \ldots, n$, \textbf{Do}\linebreak
$X'_i=X_i \cap N_i(X)$.

Output $B^*=X'_1 \times X'_2 \times \cdots \times X'_n$.
\end{algorithm}

For interval matrix multiplication in Algorithm~\ref{alg3}, it is different to common matrix multiplication since all elements in interval matrices are interval. We present the detail of interval matrix multiplication with Algorithm~\ref{alg4}.

\begin{algorithm}[H]
\caption{Interval matrix multiplication}
\label{alg4}
\KwIn{two interval matrices $A\in \mathbb{IR}^{m \times k}, B \in \mathbb{IR}^{k \times n}$, where each elements is an interval like $[\underline{d_{ij}}, \overline{d_{ij}}]$.}
\KwOut{the interval matrix $C=A \times B$.}
Initialization interval matrix $C$.

\textbf{For} {$i=1, 2, \ldots, m$}\\
\ \ \textbf{For} {$j=1, 2, \ldots, n$} \\
\ \ \ \ \textbf{For} $u=1, 2, \ldots, k$\linebreak
$$\underline{C_{ij}} = \underline{C_{ij}}+\min\{\underline{A_{iu}} \overline{B_{uj}},\ \underline{A_{iu}}  \underline{B_{uj}}, \ \overline{A_{iu}}  \underline{B_{uj}},\ \overline{A_{iu}}  \overline{B_{uj}} \},$$
$$\overline{C_{ij}}= \overline{C_{ij}}+\max\{\underline{A_{iu}}  \overline{B_{uj}},\ \underline{A_{iu}}  \underline{B_{uj}},\ \overline{A_{iu}}  \underline{B_{uj}},\ \overline{A_{iu}} \overline{B_{uj}} \}.$$\

Output $C$.
\end{algorithm}

\subsection{Backtracking algorithm based on GPGPU}
Through several successive rounds of branch and bound methods, the infeasible configurations are removed while the feasible ones are left. Since the box of these left configurations is too tiny to be recognized, we can do backtracking to get "bigger" and readable configurations. Because of bisection in the successive branch and bound method, we present the Algorithm~\ref{alg5} to demonstrate the backtracking process of feasible configurations from $n$-th round up to $(n-1)$-th round.

\begin{algorithm}[H]
\caption{Backtracking algorithm on GPGPU}
\label{alg5}
 \KwIn{A configuration of $n$-th round, $[[\underline{a_1}, \overline{a_1}], [\underline{a_2}, \overline{a_2}], \ldots, [\underline{a_n}, \overline{a_n}]] $.}
 \KwOut {A configuration of $n-1$-th round.}
\textbf{For} {each variable interval}  \textbf{Do} \linebreak
{$x_i=[\underline{a_i}, \overline{a_i}],$\linebreak
width $w_i = 2(\overline{a_i} - \underline{a_i})$,\linebreak
$x=\underline{a_i}$/$w_i$,\linebreak
$y=\overline{a_i}$/$w_i$. }

The $(n-1)$th interval $b_i$ of $x_i$ can be deduced with
 \quad $\mbox{$\underline{b_i}$}=\left\{
\begin{array}{ll}
\mbox{($\underline{a_i}+\overline{a_i}$)/2} ; & \mbox{x is integer}, \\
\underline{a_i} ;& \mbox{x is not integer}. \\
\end{array}
\right.\
$\linebreak
 \quad $ \mbox{$\overline{b_i}$}=\left\{
\begin{array}{ll}
\overline{a_i}; & \mbox{y is integer}, \\
\mbox{($\underline{a_i}+\overline{a_i})$/2} ; & \mbox{y is not integer}. \\
\end{array}
\right.\
$\\

output the new configuration $[[\underline{b_1}, \overline{b_1}], [\underline{b_2}, \overline{b_2}], \ldots, [\underline{b_n}, \overline{b_n}]]$.
\end{algorithm}

Obviously, multiple configurations of $n$-th rounds can restore up to the same configuration of $(n-1)$-th round. So we can remove the redundant configurations of $(n-1)$-th rounds and execute the backtracking algorithm successively to get more readable configurations.

\subsection{Algorithm analysis}
Our algorithm is based on branch and bound method to globally search the feasible solutions from initial space. Since the computation involves floating-point arithmetic, we use interval replacing number and use rounding down the lower bound and rounding up the upper bound, to ensure computation correctness. Hansen-Sengupta method is used to determine the root and reduce the solution space. The backtracking algorithm can isolate finite real roots. Therefore, our algorithm is right and complete. Additionally, the nonlinear algebraic equations studied in this paper have finite real solution, and the initial space is limited. So the algorithm will stop when the interval width of feasible configurations meets the desired goal. 

In branch and bound method, $n$ threads compute function interval evaluation in parallel. From view of parallel computing, the time cost of this algorithm is $O(n^2)$. Similarly, the time complexity for each thread in Hansen-Sengupta algorithm is inverse matrix computing complexity $O(n^3)$, plus with interval matrix multiplication complexity $O(n^3)$, and Gaussian elimination complexity $O(n)$. So the total time complexity of Hansen-Sengupta algorithm is $T(n)=2O(n^3)+O(n)$, as $O(n^3)$.

It is noted that our parallel algorithm does not decrease time complexity. Because we adopt coarse parallel pattern to improve whole performance. Each thread needs do whole calculation independently. There is no message communication and synchronization between different threads. This saves almost whole cost of communication.

\section{Experiment}
\label{expriment}
The detailed information of computational resource in our experiments is Intel Core i7-7700 CPU of 3.60GHz, 16GB memory, and GPGPU NVIDIA GeForce 1070 GTX with 8G video memory, 15 Multiprocessors, 1920 Cuda Cores. The version of \textit{CUDA} programming is 7.5.

To show the feasibility and effectiveness of our algorithm, we collect $55$ equations from real applications and academic references with one constraint that the variable count is not more than $9$.

We select several typical equations for analysis and detailed description. The full equation descriptions and statistics are showed in Table \ref{tab9},\ref{tab10}.

\subsection{ Analysis for typical equations }
\textbf{Example 1} Katsura8, a 9-dimensional economics problem, consider the system F(x1, \ldots, x9):\
\begin{small}
\begin{eqnarray}
\left\{ \begin{aligned}
2x_9+2x_8+2x_7+2x_6+2x_5+2x_4+2x_3+2x_2+2x_1-1=0, \\
2x_9^2+2x_8^2+2x_7^2+2x_6^2+2x_5^2+2x_4^2+2x_3^2+2x_2^2+x_1^2-x_1=0, \\
2x_9x_8+2x_8x_7+2x_7x_6+2x_6x_5+2x_5x_4+2x_4x_3+2x_3x_2+2x_2x_1-x_2=0, \\
2x_9x_7+2x_8x_6+2x_7x_5+2x_6x_4+2x_5x_3+2x_4x_2+2x_3x_1+x_2^2-x_3=0, \\
2x_9x_6+2x_8x_5+2x_7x_4+2x_6x_3+2x_5x_2+2x_4x_1+2x_3x_2-x_4=0, \\
2x_9x_5+2x_8x_4+2x_7x_3+2x_6x_2+2x_5x_1+2x_4x_2+x_3^2-x_5=0, \\
2x_9x_4+2x_8x_3+2x_7x_2+2x_6x_1+2x_5x_2+2x_4x_3-x_6=0, \\
2x_9x_3+2x_8x_2+2x_7x_1+2x_6x_2+2x_5x_3+x_4^2-x_7=0, \\
2x_9x_2+2x_8x_1+2x_7x_2+2x_6x_3+2x_5x_4-x_8=0. \\
\end{aligned} \right.
\end{eqnarray}
\end{small}

The initial interval for each variable is $[-1, 1]$.\
Using the branch and bound method iteratively, the number of interval boxes is decreasing after eight rounds. The statistics are showed in Table \ref{tab1}:\
\begin{table}[H]
\caption{B\&B result of katsura8}
\label{tab1}
\begin{tabular*}{8.5cm}{llll}
\hline\noalign{\smallskip}
round &  boxes & time(s) &width of interval   \\
\hline\noalign{\smallskip}
1 & 1 & 0 & 2 \\
2 & 511 & 0.001 & 1 \\
3 & 4847 & 0.01 & 0.5 \\
4 & 92880 & 0.2 & 0.25 \\
5 & 2549461 & 8 & 0.125 \\
6 & 76736464 & 38 & 0.0625 \\
7 & 794760638 & 6556 & 0.03125 \\
8 & 430650061 & 11203 & 0.015625\\
\noalign{\smallskip}\hline
\end{tabular*}
\end{table}

We continue to do branch and bound method and use Hansen-Sengupta method after the 8th round. We define the boxes number after Hansen-Sengupta method as new boxes. The statistics are showed in Table \ref{tab2}:

\begin{table}[H]
\caption{Hansen-Sengupta result of katsura8}
\label{tab2}
\begin{tabular*}{8.5cm}{llll}
\hline\noalign{\smallskip}
round & origin boxes  &  new boxes  &time(s)  \\
\hline\noalign{\smallskip}
8 & 430650061 & 4900601 & 1271\\
9 & 30098623 & 4347544 & 137\\
10 & 42617540 & 5532835 & 310\\
11 & 64267381 & 7946423 & 380\\
12 & 105056246  & 11571157 & 478\\
13 & 165616092 & 15474734 & 750\\
14 & 243760141 & 18711270 & 1013\\
15 & 295102270 & 20306174 & 1359\\
16 & 343840882 & 20242072 & 1624\\
\noalign{\smallskip}\hline
\end{tabular*}
\end{table}

The number of interval boxes is basically stable, then we can merge the interval boxes by backtracking algorithm, and isolate $36$ real interval roots, width is $0.125$. The statistics are showed in Table \ref{tab3},

\begin{table}[H]

\caption{backtracking result of katsura8}
\label{tab3}
\begin{tabular*}{8.5cm}{lll}
\hline\noalign{\smallskip}
 backtracking boxes  & time(s) &width of interval   \\
\hline\noalign{\smallskip}
20242072  & 29 & 0.0006103515635\\
18711270  & 26 & 0.001220703135\\
11511157 & 18 & 0.00244140625\\
7665686  & 13 & 0.0048828125\\
2221856  & 7 & 0.00976525015\\
456692  & 2 & 0.01953125\\
76643  & 0.3 & 0.00390625\\
14006  & 0.01 & 0.00390625\\
2957  & 0.001 & 0.0078125\\
691  & 0.001 & 0.015625\\
246  & 0.001 & 0.03125\\
79  & 0.001 & 0.0625\\
36  & 0.001 & 0.125\\
36  & 0.001 & 0.25\\
\noalign{\smallskip}\hline
\end{tabular*}
\end{table}

\textbf{Example 2}  Noon9, a neural network Lotka-Volterra system. Consider the system F(x1, \ldots, x9)\
\begin{small}
\begin{equation}
 \left\{ \begin{aligned}
x_1x_2^2+x_1x_3^2+x_1x_4^2+x_1x_5^2+x_1x_6^2+x_1x_7^2+x_1x_8^2+x_1x_9^2-1.1x_1+1=0, \\
x_2x_1^2+x_2x_3^2+x_2x_4^2+x_2x_5^2+x_2x_6^2+x_2x_7^2+x_2x_8^2+x_2x_9^2-1.1x_2+1=0, \\
x_3x_1^2+x_3x_2^2+x_3x_4^2+x_3x_5^2+x_3x_6^2+x_3x_7^2+x_3x_8^2+x_3x_9^2-1.1x_3+1=0, \\
x_4x_1^2+x_4x_2^2+x_4x_3^2+x_4x_5^2+x_4x_6^2+x_4x_7^2+x_4x_8^2+x_4x_9^2-1.1x_4+1=0, \\
x_5x_1^2+x_5x_2^2+x_5x_3^2+x_5x_4^2+x_5x_6^2+x_5x_7^2+x_5x_8^2+x_5x_9^2-1.1x_5+1=0, \\
x_6x_1^2+x_6x_2^2+x_6x_3^2+x_6x_4^2+x_6x_5^2+x_6x_7^2+x_6x_8^2+x_6x_9^2-1.1x_6+1=0, \\
x_7x_1^2+x_7x_2^2+x_7x_3^2+x_7x_4^2+x_7x_5^2+x_7x_6^2+x_7x_8^2+x_7x_9^2-1.1x_7+1=0, \\
x_8x_1^2+x_8x_2^2+x_8x_3^2+x_8x_4^2+x_8x_5^2+x_8x_6^2+x_8x_7^2+x_8x_9^2-1.1x_8+1=0, \\
x_9x_1^2+x_9x_2^2+x_9x_3^2+x_9x_4^2+x_9x_5^2+x_9x_6^2+x_9x_7^2+x_9x_8^2-1.1x_9+1=0. \\
\end{aligned} \right.
\end{equation}
\end{small}

The initial interval for each variable is $[-8, 8]$.
Using the branch and bound method iteratively, the number of interval boxes is decreasing after ten rounds.

The statistics are showed in Table \ref{tab3}:\
\begin{table}[H]
\caption{B\&B result of noon9}
\label{tab4}
\def\tabblank{\hspace*{10mm}} 
\begin{tabular*}{8.5cm}{llll}
\hline\noalign{\smallskip}
round & boxes & time(s) &width of interval   \\
\hline\noalign{\smallskip}
1 & 1 & 0 & 16 \\
2 & 512 & 0.001 & 8 \\
3 & 5120 & 0.001 & 4 \\
4 & 42176 & 0.01 & 2 \\
5 & 13511 & 0.2 & 1 \\
6 & 222941 & 1.5 & 0.5 \\
7 & 1156044 & 8 & 0.25 \\
8 & 11471246 & 102 & 0.125\\
9 & 33925040 & 300 & 0.0625\\
10 & 28775053 & 285 & 0.03125\\
\noalign{\smallskip}\hline
\end{tabular*}
\end{table}

We continue to do branch and bound method and use Hansen-Sengupta method after the 10th round. The statistics are showed in Table \ref{tab5}:

\begin{table}[H]
\caption{ Hansen-Sengupta result of noon9}
\label{tab5}
\begin{tabular*}{8.5cm}{llll}
\hline\noalign{\smallskip}
round &  origin boxes  & new boxes  &time(s)  \\
\hline\noalign{\smallskip}
10 & 28775053 & 477 & 78\\
\noalign{\smallskip}\hline

\end{tabular*}
\end{table}

It is easily seen that the Hansen-Sengupta method is efficient for the function., then we can merge the interval boxes by backtracking algorithm, and isolate $19$ real interval roots, width is $0.0625$. The statistics are showed in Table \ref{tab6}.

\begin{table}[H]
\caption{backtracking result of noon9}
\label{tab6}
\begin{tabular*}{8.5cm}{lll}
\hline\noalign{\smallskip}
 backtracking boxes & time(s) &width of interval   \\
\hline\noalign{\smallskip}
477 & 0.001 & 0.03125\\
19 & 0.001 & 0.0625\\
19 & 0.001 & 0.125\\
\noalign{\smallskip}\hline
\end{tabular*}
\end{table}

\textbf{Example 3}  Kinema, a robot kinematics problem. Consider the system $F(x1, \ldots, x9)$:\
\begin{small}
\begin{eqnarray}
\left\{ \begin{aligned}
x_1^2+x_2^2+x_3^2-12x_1-68=0, \\
x_4^2+x_5^2+x_6^2-12x_5-68=0, \\
x_7^2+x_8^2+x_9^2-24x_8-12x_9+100=0, \\
x_1x_4+x_2x_5+x_3x_6-6x_1-6x_5-52=0, \\
x_1x_7+x_2x_8+x_3x_9-6x_1-12x_8-6x_9+64=0, \\
x_4x_7+x_5x_8+x_6x_9-6x_5-12x_8-6x_9+32=0, \\
2x_2+2x_3-x_4-x_5-2x_6-x_7-x_9+18=0, \\
x_1+x_2+2x_3+2x_4+2x_6-2x_7+x_8-x_9-38=0, \\
x_1+x_3-2x_4+x_5-x_6+2x_7-2x_8+8=0. \\
\end{aligned} \right.
\end{eqnarray}
\end{small}

The initial interval for each variable is $[-32, 32]$.
Using the branch and bound method iteratively, the number of interval boxes is decreasing after ten rounds. The statistics are showed in Table \ref{tab7}:\
\begin{table}[H]
\caption{B\&B result of kinema}
\label{tab7}
\begin{tabular*}{8.5cm}{llll}
\hline\noalign{\smallskip}
round &  boxes  & time(s) &width of interval   \\
\hline\noalign{\smallskip}
1 & 1 & 0 & 64 \\
2 & 288 & 0.001 & 32 \\
3 & 5714 & 0.01 & 16 \\
4 & 56036 & 0.3 & 8 \\
5 & 407380 & 3 & 4 \\
6 & 4369735 & 28 & 2 \\
7 & 26247657 & 243 & 1 \\
8 & 87424465 & 1284 & 0.5\\
9 & 33925040 & 300 & 0.25\\
10 & 28775053 & 285 & 0.125\\
11 & 238223456 & 14930 & 0.0625\\
12 & 187213894 & 12250 & 0.03125\\
13 & 158874057 & 8450 & 0.015625\\
14 & 155682210 & 7321 & 0.0078125\\
15 & 153442353 & 7022 & 0.00390625\\
16 & 157313456 & 7127 & 0.001953125\\
\noalign{\smallskip}\hline

\end{tabular*}
\end{table}

We can merge the interval boxes by backtracking algorithm, and isolate $8$ real interval roots width is $0.03125$. The statistics are showed in Table \ref{tab8}:\
\begin{table}[H]
\footnotesize
\caption{backtracking result of kinema}
\label{tab8}
\begin{tabular*}{8.5cm}{lll}
\hline\noalign{\smallskip}
  backtracking boxes  & time(s) &width of interval   \\
\noalign{\smallskip}\hline\noalign{\smallskip}
6587702 & 24 & 0.001953125\\
328906 & 1.5 & 0.00390625\\
14876 & 0.1 & 0.0078125\\
453 & 0.01 & 0.015625\\
8 & 0.001 & 0.03125\\
8 & 0.001 & 0.0625\\
\noalign{\smallskip}\hline
\end{tabular*}
\end{table}

In Table \ref{tab9}, we present the description of whole test set. In Table \ref{tab10}, we list the variable dimension, iteration rounds, number of isolated solutions of backtracking algorithm, and cost time. Some equations cannot be computed in finitely reasonable time and space, and we mark them as blank. See the full statistics for all equations at: http://github.com/eviloan/Nonlinear-Equation-Test-Result.

\begin{table}[H]
\centering
\footnotesize
\caption{description of test set}
\label{tab9}
\begin{tabular}{lll}
\hline\noalign{\smallskip}
name  & variable  &  description   \\
\noalign{\smallskip}\hline\noalign{\smallskip}
mickey & 2 & mickey-mouse example to illustrate homotopy continuation  \\
barry & 3 & barry from PoSSo  \\
arnborg & 3 & lazard, auxiliary in cyclic 7 roots, from PoSSo  \\
noon3 & 3 & neural network, Lotka-Volterra system, n=3  \\
rediff3 & 3 & 3-dimensional reaction-diffusion problem  \\
conform1 & 3 & conformal analysis of cyclic molecules, instance 1  \\
morgan & 3 & from PoSSo   \\
guakwa2 & 4 & Gaussian quadrature formula 2 knots,2 weights  \\
katsura3 & 4 & a problem of magnetism in physics n=4  \\
liu & 4 & from PoSSo   \\
cyclic4 & 4 & cyclic 4-roots problem  \\
caprasse & 4 & the system caprasse of the PoSSo test suite\\
lorentz & 4 & equilibrium of 4-dimensional Lorentz attractor  \\
moeller5 & 4 & moeller example 5, from PoSSo  \\
noon4 & 4 & neural network, Lotka-Volterra system, n=4  \\
reimer4 & 4 & the 4-dimensional system of Reimer \\
solotarev & 4 & from PoSSo  \\
katsura4 & 5 & a problem of magnetism in physics n=5  \\
noon5 & 5 & neural network, Lotka-Volterra system, n=5\\
eco5 & 5 & 5-dimensional economics problem  \\
redeco5 & 5 & reduced 5-dimensional economics problem \\
sparse5 & 5 & 5-dimensional sparse symmetric polynomial system  \\
wright & 5 & system of A.H.Wright  \\
cyclic5 & 5 & cyclic 5-roots problem  \\
reimer5 & 5 & the 5-dimensional system of Reimer \\
boon & 6 & neurophysiology, posted by Sjirk Boon  \\
eco6 & 6 & 6-dimensional economics problem  \\
katsura5 & 6 & a problem of magnetism in physics n=7  \\
noon6 & 6 & neural network, Lotka-Volterra system, n=6  \\
gaukwa3 & 6 & Gaussian quadrature formula 2 knots,2 weights  \\
trinks1 & 6 & system of Trinks from the PoSSo test suite \\
cyclic6 & 6 & cyclic 6-roots problem  \\
extcyc6  & 6 & extended cyclic 6-roots problem, to exploit the symmetry  \\
redeco6 & 6 & reduced 6-dimensional economics problem \\
romin & 6 & romin robot inverse model, from PoSSo  \\
eco7 & 7 & 7-dimensional economics problem  \\
cyclic7 & 7 & cyclic 7-roots problem  \\
redeco7 & 7 & reduced 7-dimensional economics problem \\
katsura6 & 7 & a problem of magnetism in physics n=7  \\
noon7 & 7 & neural network, Lotka-Volterra system, n=7  \\
eco8 & 8 & 8-dimensional economics problem  \\
redeco8& 8 & hand position and orientation of PUMA robot  \\
cyclic8 & 8 & cyclic 8-roots problem  \\
redeco8 & 8 & reduced 8-dimensional economics problem \\
guakwa4 & 8 & Gaussian quadrature formula 2 knots,2 weights  \\
katsura7 & 8 & a problem of magnetism in physics n=8  \\
noon8 & 8 & neural network, Lotka-Volterra system, n=8  \\
heart & 8 & heart-dipole problem  \\
SHEPWM & 8 & selective harmonic eliminated pulse width modulation \\
$s9_1$  & 8 & small system from constructive Galois theory, called $s9_1$ \\
eco9 & 9 & 9-dimensional economics problem  \\
katsura8 & 9 & a problem of magnetism in physics n=9  \\
kinema & 9 & robot kinematics problem  \\
noon9 & 9 & neural network, Lotka-Volterra system, n=9  \\
cyclic9 & 9 & cyclic 9-roots problem  \\
\noalign{\smallskip}\hline

\end{tabular}
\end{table}

\begin{table}[H]
\centering
\footnotesize
\caption{statistics of test set}
\label{tab10}
\begin{tabular}{llllll}
\hline\noalign{\smallskip}
name &  dimension &  initial interval & iterations & real roots & time(s)\\
\hline\noalign{\smallskip}
mickey & 2 & [-2 2] & 7 & 2 & 0.001  \\
barry & 3 & [-32 32] & 14 & 2 & 0.001  \\
arnborg & 3 & [-16 16] & 13 & 8 & 0.020  \\
noon3 & 3 & [-8 8] & 12 & 7 & 0.001  \\
rediff3 & 3 & [-1 1] & 9 & 2 & 0.001  \\
conform1 & 3 & [-2 2] & 3 & 0 & 0.001  \\
morgan & 3 & [-32 32] & 13 & 2 & 0.001  \\
guakwa2 & 4 & [-4 4] & 10 & 2 & 0.001  \\
katsura3 & 4 & [-1 1] & 9 & 2 & 0.001 \\
liu & 4 & [-1 1] & 9 & 1 & 0.001  \\
cyclic4 & 4 & [-16 16] & 20 & 8 & 0.430 \\
caprasse & 4 & [-4 4] & 11 & 18 & 0.020  \\
lorentz & 4 & [-2 2] & 8 & 3 & 0.001  \\
moeller4 & 4 & [-2 2] & 11 & 8 & 0.001 \\
moeller5 & 4 & [-2 2] & 11 & 6 & 0.02  \\
noon4 & 4 &  [-8 8] & 15 & 15 & 0.020  \\
reimer4 & 4 & [-1 1] & 13 &  & 0.300  \\
solotarev & 4 & [-8 8] & 14 & 6& 0.001 \\
katsura4 & 5 & [-1 1] & 8 & 8 & 0.001 \\
noon5 & 5 & [-8 8] & 13 & 11 & 9.8  \\
eco5 & 5 & [-8 8] & 8 & 4 & 0.020  \\
redeco5 & 5 & [-8 8] & 12 & 4 & 0.001  \\
sparse5 & 5 & [-2 2] & 3 & 0 & 0.001  \\
wright & 5 & [-8 8] & 8 & 32 & 0.030  \\
cyclic5 & 5 & [-16 16] & 11 & 10 & 0.05  \\
reimer5 & 5 & [-1 1] &  &  &   \\
boon & 6 & [-2 2] & 9 & 8 & 0.01  \\
eco6 & 6 & [-8 8] & 9 & 4 & 0.004  \\
katsura5 & 6 & [-1 1] & 8 & 8 & 0.03  \\
noon6 & 6 & [-8 8] & 13 & 13 & 120.6  \\
gaukwa3 & 6 & [-4 4] & 12 & 0 & 66  \\
trinks1 & 6 & [-8 8] & 11 & 8 & 0.003  \\
cyclic6 & 6 & [-16 16] & 15 & 24 & 42 \\
extcyc6  & 6 & [-16 16] &  &  &   \\
redeco6 & 6 & [-16 16] & 12 & 4 & 0.02  \\
romin & 6 & [-1 1] & 13 & 4 & 0.003  \\
eco7 & 7 & [-8 8] & 8 & 8 & 0.02  \\
cyclic7 & 7 & [-16 16] & 15 &  & 3722  \\
redeco7 & 7 & [-16 16] &  &  &   \\
katsura6 & 7 & [-1 1] & 8 & 16 & 36.2 \\
noon7 & 7 & [-8 8] & 12 & 15 & 264  \\
eco8 & 8 & [-8 8] & 7 & 8 & 0.05  \\
puma & 8 & [-1 1] & 13 & 16 & 0.02  \\
cyclic8 & 8 & [-16 16] &  &  &   \\
redeco8 & 8 & [-16 16] &  &  &   \\
guakwa4 & 8 & [-16 16] &  &  &   \\
katsura7 & 8 & [-1 1] & 9 & 16 & 7150  \\
noon8 & 8 & [-8 8] & 11 & 17 & 924  \\
heart & 8 & [-16 16] &  &  &  \\
SHEPWM & 8 & [-1 1] & 5 &  &  \\
$s9_1$ & 8 & [-8 8] & 11 & 4 &0.001  \\
eco9 & 9 & [-8 8] &  &  &   \\
katsura8 & 9 & [-1 1] & 16 &  & 13360  \\
kinema & 9 & [-32 32] & 16& 8 & 56334 \\
noon9 & 9 & [-8 8] & 10 & 19 & 1021  \\
cyclic9 & 9 & [-16 16] &  &  &   \\
\noalign{\smallskip}\hline

\end{tabular}
\end{table}

\section{Conclusions}
\label{final}
In this paper, we propose a global interval algorithm, to solve the nonlinear algebraic equations and get the feasible interval solutions. Our algorithm is implemented on GPGPU. Our future work is to improve the computational efficiency of the algorithm. One of the hopeful direction is finding an optimization method to quickly find the maximum and the minimum of the function in a given region to obtain a more tighter interval, improve the speed of branch and bound method.




\end{document}